\documentclass[12pt,twoside,reqno]{article}
\usepackage{amsfonts}
\usepackage{graphicx}

\def\bkR{{\rm I\kern-.17em R}}
\def\sq{{\rm \sqcap \kern-.62em _{-} }}
\begin{document}
\begin{center}
{{\large \bf A Sparse-Sparse Iteration for Computing a Sparse
Incomplete Factorization of the Inverse of an SPD Matrix
 }} \vspace*{0.5cm}
\\
{{\bf Davod Khojasteh Salkuyeh}$^1$ {\bf and} {\bf Faezeh
Toutounian}$^2$} \vspace*{0.3cm}
\\
$^1$ {\small Department of Mathematics, University of Mohaghegh Ardabili,\\
P. O. Box. 56199-11367, Ardabil, Iran\\
E-mail: khojaste@uma.ac.ir \\}

$^2$ {\small School of  Mathematical Sciences, Ferdowsi University of Mashhad,\\
P. O. Box. 1159-91775, Mashhad, Iran  \\
E-mail: toutouni@math.um.ac.ir \\}

\end{center}

\bigskip

\medskip

\noindent {\bf Abstract:} In this paper, a method via
sparse-sparse iteration for computing a sparse incomplete
factorization of the inverse of a symmetric positive definite
matrix is proposed. The resulting factorized sparse approximate
inverse is used as a preconditioner for solving symmetric positive
definite linear systems of equations by using
the preconditioned conjugate gradient algorithm. Some numerical experiments on test
matrices from the Harwell-Boeing collection for comparing the
numerical performance of the presented method with one available
well-known algorithm are also given.

\bigskip
\noindent{\small\bf AMS Subject Classification :} 65F10, 65F50.\\
\noindent{\textit{Keywords}}: Sparse matrices; Factorized sparse
approximate inverse; Preconditioning; Krylov subspace methods;
Symmetric positive definite; Preconditioned  CG algorithm

\bigskip

\bigskip

\medskip

\noindent \textbf{1. Introduction}

\bigskip

\bigskip

Consider the nonsingular linear system of equations
\begin{equation}\label{e:eq1}
    Ax=b,\\
\end{equation}
where the coefficient matrix  $A \in \mathbb{R}^{n \times n}$ is
large, sparse and $x,b\in \mathbb{R}^{n}$. It is well-known that
the rate of convergence of iterative methods such as Krylov
subspace methods for solving (\ref{e:eq1}) is strongly influenced
by the spectral properties of $A$. Hence, iterative methods
usually involve a second matrix that transforms the coefficient
matrix into one with a more favorable spectrum. The transformation
matrix is called a preconditioner.
If $M$ is a nonsingular matrix that approximates the inverse of $A$ ($M \approx A^{-1})$, then
the transformed linear system
\begin{equation}\label{e:eq2}
    AMy=b,\quad  x=My,\\
\end{equation}
will have the same solution as system (\ref{e:eq1}), but the convergence rate of iterative methods applied to (\ref{e:eq2}) may be higher. System (\ref{e:eq2}) is preconditioned from the right, but left preconditioning is also possible, i.e., $MAx=Mb$. One can also define split-preconditioned systems.
Let us assume that $A$ has the LU factorization and
\begin{equation}\label{e:eq3}
M=M_UM_L, ~~~~{\mbox{where}}~~ M_U \approx U^{-1}~{\mbox{and}}~ M_L
\approx L^{-1},
\end{equation}
where $L$ and $U$ are the lower and upper triangular factors of
$A$. This type of preconditioning is known as factorized
approximate inverses and $M_U$ and $M_L$ are called approximate
inverse factors of $A$. Here, the transformed linear system can be
considered as follows
\begin{equation}\label{e:splsys}
    M_UAM_Ly=M_Ub,~~x=M_Ly.\\
\end{equation}
\noindent System (\ref{e:splsys}) is called a split-preconditioned system.

\indent In this paper we focus our attention on the computation of
sparse approximate inverse factors of a matrix. There are
different ways to compute sparse approximate inverse factors of a
matrix and each of them has its own advantages and disadvantages.
In \cite{SYMAINV,AINV}, the AINV method was proposed which is
based on an algorithm which computes two sets of vectors $\{ z_i
\}_{i=1}^{n}$ and  $\{ w_i \}_{i=1}^{n}$ which are
$A$-biconjugate, i.e., such that $w_i^TAz_j=0$ if and only if $i
\neq j$. Although the construction phase for the original AINV
algorithm is sequential, its application is highly parallel,
since it consists of  matrix-vector products. A fully parallel AINV
algorithm can be achieved by means of graph partitioning (see
\cite{SURVEY, Marin,Markov}). For symmetric positive definite
(SPD) matrices, there exists a variant of the AINV method,
denoted by SAINV (for Stabilized AINV), that is breakdown-free
\cite{SAINV}. Another approach which was proposed by Kolotilina
and Yeremin is the FSAI algorithm \cite{Kol1,Kol2}. They assume
that $A$ is SPD and then construct factorized sparse approximate
inverse preconditioners which are also SPD. Each factor
implicitly approximates the inverse of the lower triangular
Cholesky factor of $A$. This method can be easily  extended to
the nonsymmetric case. The FSAI algorithm is inherently parallel
but its main disadvantage is the need to prescribe the sparsity of
approximate inverse factors in advance.

\indent In this paper,  we first propose an iterative method for
solving SPD linear systems of equations, and then, by exploiting this method, we
develop an algorithm for computing an incomplete factorization of
the inverse of an SPD matrix. The resulting factorized sparse
approximate inverse is used as an explicit preconditioner for the
solution of $Ax=b$ by the preconditioned conjugate gradient (PCG)
method.

\indent Throughout the paper $\|z\|_A$ stands for the $A$-norm of
any vector $z$, i.e., $\|z\|_A=(Az,z)^{1/2}$. We will denote the
largest and smallest eigenvalues of the matrix $X$ by
$\lambda_{max}(X)$ and $\lambda_{min}(X)$, respectively.

 \indent This paper is organized as follows. In
section 2, we introduce an approach for computing a sparse
approximate solution of an SPD linear system of equations.
Section 3 is devoted to computing an incomplete factorization of the
inverse of an SPD matrix. Numerical experiments are given in
section 4. Finally, we give some concluding remarks in section 5.

\bigskip

\bigskip

\medskip

\noindent  \textbf{2. Sparse approximate solution of an SPD linear
system of equations}

\bigskip

\bigskip

\indent In this section, we first present an approach,
based on the projection method, for solving an SPD linear system of
equations. Then we develop an algorithm for computing a sparse
approximate solution of an SPD linear system of equations.

\indent Let $A=(a_{ij})$ be an SPD matrix and let us consider a
projection method with $\mathcal{L}=\mathcal{K}={\mbox{span}}
\{e_{i_1},e_{i_2},\ldots,e_{i_m}\}$, where $e_{i_j}$ is the
$i_j$-th column of identity matrix and $m$ is a small natural
number. Given an initial guess $x$ of the solution of
(\ref{e:eq1}) and the residual vector $r=b-Ax$,  the new
approximation takes the form
\begin{equation}\label{e:eq5}
   x_{new}=x+Ey,
\end{equation}
for some $y \in \mathbb{R}^m$, and
$E=[e_{i_1},e_{i_2},\ldots,e_{i_m}]$. The Petrov-Galerkin
condition $ r-AEy \perp \mathcal{L}$ yields
\begin{equation}\label{e:eq6}
y=(E^TAE)^{-1}E^Tr.
\end{equation}
As known \cite{Saadbook}, this kind of update minimizes
$$ \| x+E \widetilde{y}-x_{exact} \|_A$$
over all $\widetilde{y }\in \mathbb{R}^m$, where $x_{exact}$ is
the exact solution of $Ax=b$. It is obvious that the matrix
$S=E^T A E$ is an SPD matrix of dimension $m$. Defining
$\mathcal{J}=\{i_1,i_2,\ldots,i_m \}$, the matrix $E^TAE$ is the
principal submatrix of $A$ consisting of the rows and columns
whose indices are in $\mathcal{J}$. This new approach for solving
an SPD linear system of equations can be stated as follows.

\bigskip

\medskip

\underline{{\bf Algorithm 1:}}\\
\vspace{-0.6cm}

\begin{itemize}
\item[1.] Choose  an initial guess $x$ and compute $r=b-Ax$
\vspace{-0.2cm}

\item[2.] Until convergence, Do \vspace{-0.2cm}

\item[3.] ~~~~~~~Select $\mathcal{J}=\{i_1,i_2,\ldots,i_m \}
\subseteq \{1,2,\ldots,n \}$\vspace{-0.3cm}

\begin{itemize}
    \item[]\hspace{-0.1cm} and $E:=[e_{i_1},e_{i_2},\ldots,e_{i_m}]$ \vspace{-0.2cm}
\end{itemize}

\item[4.]~~~~~~~Solve $(E^TAE)y=E^Tr$ for $y$ \vspace{-0.2cm}

\item[5.]~~~~~~~Compute $x:=x+Ey$ \vspace{-0.2cm}

\item[6.]~~~~~~~Compute $r:=r-AEy$ \vspace{-0.2cm}

\item[7.] EndDo
\end{itemize}

\bigskip

Step 6 of this algorithm  can be written as
\begin{equation}\label{e:eq7}
 r:=r-\sum_{k \in
\mathcal{J}} y_k~a_{:,k},
\end{equation}
 where $a_{:,k}$ is the $k$-th column of
$A$ and $y=[y_{1},y_{2},\ldots,y_{m}]^T$. The relation
(\ref{e:eq7}) shows that, for updating $r$, we need $m$ sparse
SAXPY operations (a SAXPY operation is defined as $z:=x+ \alpha
y$, where $x$ and $y$ are $n$-vectors and $\alpha$ is a scalar).
The following theorem regarding the convergence rate of the
Algorithm 1 can be stated.

\bigskip

\smallskip

\noindent {\bf Theorem 1.} {\it Let $A$ be a symmetric positive
definite matrix.  Assume that, at each projection step, the
selected index set $ \{ i_1,i_2,\ldots,i_m \} $ contains the indices
of the $m$ components with largest absolute value in the current
residual vector $r=b-Ax$. Then
\begin{equation}\label{e:eq8}
\| d \|_A^2-\| d_{new} \|_A^2 \geq \frac{\sum_{k \in
\mathcal{J}}r_k^2}{\sum_{k \in \mathcal{J}}a_{kk}},
\end{equation}
and
\begin{equation}\label{e:eq9}
\| d_{new} \|_A \leq (1-(\frac{\lambda_{min}(A)}{\sum_{k \in
\mathcal{J}}a_{kk}})( \frac{\sum_{k \in
\mathcal{J}}r_{k}^2}{\sum_{k=1}^{n}r_k^2}))^{1/2}\| d
\|_A,\\
\end{equation}
\noindent where $d_{new}=A^{-1}b-x_{new}$ and $d=A^{-1}b-x$.
Relation (\ref{e:eq9}) shows that Algorithm 1 converges for any
initial guess.}

\bigskip
\noindent {\bf Proof.}  We start by observing that $d_{new}=d-Ey$,
 $Ad=r$, and
$$ (Ad_{new},d_{new}) = (Ad,d)-(y,E^Tr).$$
From  (\ref{e:eq6}) and using the Courant-Fisher min-max theorem
\cite{Axe1,Saadbook},  we have
\begin{eqnarray*}
(y,E^Tr) & = & ((E^TAE)^{-1}E^Tr,E^Tr) \\
& \geq  & \frac{\|E^Tr\|_2^2}{\lambda_{max}(E^TAE)} \\
& \geq  & \frac{\|E^Tr\|_2^2}{\sum_{k \in \mathcal{J}}a_{kk}},
\end{eqnarray*}
and \begin{eqnarray*} (Ad,d) =  (r,A^{-1}r) \leq
\frac{\|r\|_2^2}{\lambda_{min}(A)}.
\end{eqnarray*}
From these observations the desired results immediately follow.
Relation (\ref{e:eq9}) establishes the convergence of the method,
since $\sum_{k \in \mathcal{J}}r_{k}^2\leq \sum_{k=1}^{n}r_k^2$ and
$\lambda_{min}(A)\leq \lambda_{min}(E^TAE)\leq \sum_{k \in
\mathcal{J}}a_{kk}$ (see \cite{Axe1}). $\hspace{1cm}\Box$
\bigskip

\smallskip

\indent  This theorem not only shows the convergence of the
algorithm but also the rate of the reduction in the square of the
$A$-norm of the error (Eq. (\ref{e:eq8})). In fact, the indices
$i_j, j=1, \ldots , m$ are chosen in such a way that the reduction
in the square of the $A$-norm of the error is as large as possible.
If $A$ is a symmetric diagonally scaled matrix then $\sum_{k\in
\mathcal{J}}a_{kk}=m$ and in the Eqs. (\ref{e:eq8}) and
(\ref{e:eq9}), $\sum_{k\in \mathcal{J}}a_{kk}$  may be replaced by
$m$.

\indent Now, by using Algorithm 1,  we propose an algorithm to
compute a sparse approximate solution of an SPD linear system of
equations. In this algorithm no dropping strategy is needed and
sparsity of the solution is preserved only by specifying the maximum
number of its nonzero entries, $lfil$, in advance. In each iteration
at most $m$ $(m \ll n)$ entries are added to the current approximate
solution. This algorithm can be stated as follows.

\bigskip

\medskip

\underline{ {\bf Algorithm 2 :} Sparse approximate solution to the SPD system $Ax=b$ }\\
\vspace{-0.7cm}

\begin{itemize}
\item[1.] Set $x:=0$ and $r:=b$ \vspace{-0.2cm}

\item[2.] While $\|\ r \| > {eps}$ and $nnz(x) < lfil$
Do \vspace{-0.2cm}

\item[3.]~~~~~~~Select the indices of $m$ components with largest
absolute value in the \vspace{-0.3cm}

\begin{itemize}
\item[]current residual vector $r$, i.e.,
$\mathcal{J}=\{i_1,i_2,\ldots,i_m \} \subseteq \{1,2,\ldots,n
\}$\\
\hspace{0.1cm}and set $E:=[e_{i_1},e_{i_2},\ldots,e_{i_m}]$
\vspace{-0.1cm}

\end{itemize}
\vspace{-0.35cm}

\item[4.]~~~~~~~Solve $(E^TAE)y=E^Tr$ for $y$ \vspace{-0.2cm}

\item[5.]~~~~~~~Compute $x:=x+Ey$ \vspace{-0.2cm}

\item[6.]~~~~~~~Compute $r:=r-AEy$ \vspace{-0.2cm}

\item[7.] EndDo
\end{itemize}

\bigskip

The vector $x$ computed by Algorithm 2 has at most $lfil$ nonzero
entries. In practical implementations of Algorithm 2 the number $m$
is usually chosen to be too small, for example $m=1,2 ~{\mbox{or}}~3$.
Throughout this paper we take $m=2$.  The parameter $eps$ is used
for stopping the process when the residual norm is small enough.
As can be seen, in Algorithm 2 no dropping strategy is used and
in each step of the algorithm, according to Theorem 1, the
$A$-norm of the error is reduced.

\bigskip

\bigskip

\medskip

\noindent \textbf{3. Approximate inverse factors of a matrix via
sparse-sparse iterations}

\bigskip

\bigskip

In this section, for computing a sparse factorized approximate
inverse of an SPD matrix, we combine Algorithm 2 of section 2 with
the AIB (Approximate Inverse via Bordering) algorithm proposed by
Saad in \cite{Saadbook}. We first give a brief description of
the AIB algorithm for symmetric matrices.

In the AIB algorithm, the sequence of matrices
\begin{equation}\label{e:eqs6}
A_{k+1}=\left ( \begin{array}{cc}
A_k & v_k\\
v_k^T & \alpha_{k+1} \\
\end{array} \right),
\end{equation}
is made in which $A_n=A$. If the inverse factor $U_k$ is available
for $A_k$, i.e.,
\begin{equation}\label{e:eqs6p}
U_k^TA_kU_k=D_k,
\end{equation}
then the inverse factor $U_{k+1}$ for $A_{k+1}$ will be obtained
by writing
\begin{equation}\label{e:eqs7}
\left ( \begin{array}{cc}
 U_k^T & 0\\
-z_k^T & 1 \\
\end{array} \right )
\left ( \begin{array}{cc}
A_k & v_k \\
v_k^T & \alpha_{k+1} \\
\end{array} \right )
\left ( \begin{array}{cc}
U_k & -z_k\\
0 & 1 \\
\end{array} \right )
=\left ( \begin{array}{cc}
D_k & 0\\
0 & \delta_{k+1} \\
\end{array} \right ),
\end{equation}
in which
\begin{equation}\label{e:eqs8}
\hspace{-1.7cm}~~~A_kz_k=v_k,
\end{equation}
\begin{equation}\label{e:eqs10}
\hspace{0.3cm}\delta_{k+1}=\alpha_{k+1}-z_k^Tv_k.
\end{equation}
Relation (\ref{e:eqs10}) can be exploited if the system
(\ref{e:eqs8}) is solved exactly. Otherwise we should use
\begin{eqnarray}\label{e:eqs11}
\nonumber  \delta_{k+1} &=& \alpha_{k+1}-v_k^Tz_k-z_k^T(v_k-A_kz_k) \\
\nonumber               &=& \alpha_{k+1}-v_k^Tz_k-z_k^Tr_k \\
                        &=& \alpha_{k+1}-z_k^T(v_k+r_k),
\end{eqnarray}

\noindent instead of (\ref{e:eqs10}), where $r_k=v_k-A_kz_k$.
Starting from $k=1$, this procedure suggests an algorithm for
computing the inverse factors of $A$. If a sparse approximate
solution of (\ref{e:eqs8}) is computed, then an approximate
factorization of $A^{-1}$ is obtained. To do this, we use
Algorithm 2 of section 2. This scheme can be summarized as
follows.

\bigskip

\indent \underline{{\bf Algorithm 3. AIB algorithm}}\\
\vspace{-0.7cm}

\begin{itemize}
\item[1.] Set $A_1=[a_{11}]$, $U_1=[1]$  and
$\delta_1=a_{11}$\vspace{-0.2cm}

\item[2.] For $k=1,\ldots,n-1$ Do: (in parallel) \vspace{-0.2cm}

\item[3.] ~~~~~~~Compute a sparse approximate solution to
$A_kz_k=v_k$,  \vspace{-0.2cm}
\begin{itemize}
\item[]\hspace{0.1cm}by using Algorithm 2,  and the residual
$r_k=v_k-A_kz_k$. \vspace{-0.2cm}
\end{itemize}

\item[4.]~~~~~~~Compute
$\delta_{k+1}=\alpha_{k+1}-z_k^T(v_k+r_k)$. \vspace{-0.2cm}

\item[5.]~~~~~~~Form $U_{k+1}$ and $D_{k+1}$\vspace{-0.2cm}

\item[6.] EndDo. \vspace{-0.2cm}

\item[7.] Set $U:=U_n$ and $D:=D_n$
\end{itemize}

\bigskip

\noindent This algorithm returns $U$ and $D$ such that $U^TAU
\approx D$. The following theorem  shows that $\delta_{k+1}$ is
always positive, independently of the accuracy with which the
system (\ref{e:eqs8}) is solved.

\bigskip

\noindent {\bf Theorem 2.} {\it Let $A$ be an SPD matrix. Then,
the scalar $\delta_{k+1}$ computed in step 4 of Algorithm 3 is
positive.}

\bigskip

\noindent {\bf Proof.}  Let $r_k$ be the residual obtained in step
3 of the AIB algorithm, i.e.,
$$r_k=v_k-A_kz_k.$$
\noindent Hence, we have
$$z_k=A_k^{-1}(v_k-r_k).$$
\noindent By a little computation one can see that
$$\delta_{k+1}=\alpha_{k+1}-v_k^TA_k^{-1}v_k+r_k^TA_k^{-1}r_k=s+r_k^TA_k^{-1}r_k,$$
\noindent where $s=\alpha_{k+1}-v_k^TA_k^{-1}v_k \in \mathbb{R}$
is the Schur complement of $A_{k+1}$ and is a positive real
number (see Theorem 3.9 in \cite{Axe1}). So, the scalar
$\delta_{k+1}$ is positive, since $A$ is an SPD matrix and
$r_k^TA_k^{-1}r_k>0$ for $r_k\not= 0$. $\hspace{1cm}\Box$

Hence the AIB algorithm is well-defined for SPD  matrices.

\bigskip

\medskip

\noindent \textbf{ 4. Numerical examples}

\bigskip

\indent All the numerical experiments presented in this section
were computed in double precision using Fortran PowerStation
version 4.0 on a Pentium 4 PC, with a 3.06 GHz CPU and
1.00GB of RAM.

\indent  For the first set of the numerical experiments, we used
nine SPD matrices (BCSSTK* and S*RMT3M*) from the Matrix-Market
website \cite{MM} and three matrices (EX15, MSC04515 and KUU )
from Tim Davis's collection \cite{Davis}. These matrices with
their generic properties are given in Table 1. For each matrix,
the problem size $n$ and the number of nonzero entries in the
lower triangular part $nnz$ are provided. In last two columns, the
number of iterations (iters) and time required to solve the
linear system using the conjugate gradient method without any
scaling are given. The time was measured with the function
\verb"etime()"  and given in seconds. The stopping criterion
\[
\frac{\| b-Ax_i \|_2}{\| b \|_2} < 10^{-8},
\]
was used and the initial guess was taken to be the zero vector. For all
the examples, the right hand side of each system was taken such that
the exact solution is a vector with random entries uniformly
distributed in $(0,1)$. No significant differences were observed for
other choices of the right hand side vector.  The maximum number of
iterations was 10000. In all the tables a dagger ($\dag$) indicates
no convergence of the iterative method.

%************Table 1  *******************

\begin{table}

\caption{ First set of test problems information.}\vspace{0.4cm}

\centering
\begin{tabular}{|l|c|c|c|c|}
\hline
~~~~matrix~~~~ & ~~~$n$ ~~~ & ~~~$nnz$~~~ & ~~ time~~ &~~ iters ~~ \\
\hline \hline

BCSSTK11   & 1473 & 17857  & - & $\dag$
\\\hline

BCSSTK13   & 2003 & 42943  & - & $\dag$
\\\hline

BCSSTK15   & 3948 & 60882  & 29.07 & 9219
\\\hline

BCSSTK21   & 3600 & 15100  & 15.32  & 7805
\\\hline

BCSSTK38   & 8032 & 181746 & - & $\dag$
\\\hline

S1RMT3M1   & 5489 & 112505 & 24.13  & 4953
\\\hline

S2RMT3M1   & 5489 & 112505 & -      & $\dag$
\\\hline

S3RMT3M1   & 5489 & 112505 & -   & $\dag$
\\\hline

S3RMT3M3   & 5357 & 106526 & -      &  $\dag$
\\\hline

MSC04515   & 4515 & 51111  & 15.12 & 4728
\\ \hline

EX15       & 6867 & 52769  & 6.48   & 1506
\\ \hline

KUU        & 7102 & 173651 & 3.84   &  550\\ \hline
\end{tabular}
\end{table}

\indent We compare the numerical results of the new preconditioner
with that of the SAINV preconditioner. The AINV and the SAINV
algorithms have been widely compared with other preconditioning
techniques, showing that they are the most effective algorithms
for computing a sparse incomplete factorization of the inverse of
a matrix \cite{SURVEY,SAINV,SYMAINV,Comp,AINV}. For the SAINV
algorithm we used the SAINV code of the SPARSLAB software
provided by Tuma \footnote{http://www.cs.cas.cz/$\sim$tuma/sparslab.html} with
drop tolerance $\tau=0.1$. This drop tolerance is very often the
right one based on the numerical results reported in several
papers. For Algorithm 2, we used the parameters $eps=0.01$ and
$lfil=10$. We also used a parameter $lfil$ such that the number
of nonzero entries in the incomplete $U$ factor divided by the
number of nonzero entries in the upper triangular part of $A$,
$\rho$, is approximately equal to or less than that of the SAINV
preconditioner. The results of the split-preconditioned CG
algorithm \cite{Saadbook} in conjunction with the SAINV
preconditioner and Algorithm 3 are given in Table 2. This table
reports the density ($\rho$), the number of split-preconditioned
CG iterations for convergence (P-Its), the setup time for the
preconditioner (P-time), the time for the split-preconditioned
iterations (It-time), and T-time which is equal to the sum of
P-time and It-time. Numerical results presented in this table
show that both algorithms are robust and the new method is better
than the SAINV algorithm for 9 out of 12 problems, especially on
the shell problems (S3RMT3M1 and S3RMT3M3). The results of this
table also indicate that the parameters $eps=0.01$ and $lfil=10$
give good results.

%************Table 2  *******************
\begin{table}
 \caption{Setup time to compute sparse approximate
inverse factors and results for the split-preconditioned CG
algorithm.}\vspace{0.4cm}

{\tiny

\centering
\begin{tabular}{|l||c|c|c|c|c|c||c|c|c|c|c|  }
\hline &\multicolumn{6}{|c||}{Algorithm 3} &
\multicolumn{5}{|c|}{SAINV Algorithm } \\
\hline \hline

~~~~matrix & $lfil$&$\rho$ & P-Its & P-time & It-time& T-time &
$\rho$&
P-Its & P-time & It-time & T-time   \\
\hline \hline

BCSSTK11 &13 & 0.58 & 628 & 0.23 & 1.20 & 1.43 & 0.58 & 2099 & 0.11 & 3.95 & 4.06\\
         &10 & 0.45 & 650 & 0.17 & 1.19 & 1.36 & & & & & \\ \hline

BCSSTK13 &29 & 0.71 & 343 & 0.65 & 1.42 & 2.07 & 0.72 & 370  & 0.63 & 1.55 & 2.18\\
         &10 & 0.26 & 550 & 0.10 & 1.80 & 1.90 & & & & & \\ \hline

BCSSTK15 &9  & 0.35 & 504 & 0.19 & 2.67 & 2.86 & 0.32 & 214  & 0.22 & 1.11 & 1.33 \\
         &10 & 0.37 & 491 & 0.20 & 2.73 & 2.93 & & & & & \\ \hline

BCSSTK21 &6  & 0.97 & 246 & 0.09 & 0.73 & 0.82 & 1.05 & 167  & 0.06 & 0.5 & 0.56\\
         &10 & 1.48 & 164 & 0.13 & 0.52 & 0.65 & & & & & \\ \hline

BCSSTK38 &11 & 0.28 & 559 & 0.70 & 7.72 & 8.42 & 0.29 & 1131 & 0.81 & 15.67 & 16.48\\
         &10 & 0.25 & 572 & 0.64 & 7.77 & 8.41 & & & & & \\ \hline

S1RMT3M1 &15 & 0.41 & 244 & 0.45 & 2.28 & 2.73 & 0.83 & 257  & 0.98 & 2.97 & 3.95\\
         &10 & 0.29 & 318 & 0.32 & 2.77 & 3.09 & & & & & \\ \hline

S2RMT3M1 &15 & 0.41 & 526 & 0.58 & 4.97 & 5.55 & 1.29 & 539  & 1.53 & 7.5  & 9.03\\
         &10 & 0.28 & 585 & 0.34 & 5.11 & 5.45 & & & & & \\ \hline

S3RMT3M1 &15 & 0.36 & 1298 & 0.65 & 11.89 & 12.54 & 2.81 & 5434 & 6.33 & 117.80 & 124.13 \\
         &10 & 0.26 & 1458 & 0.33 & 12.55 & 12.88 & & & & & \\ \hline

S3RMT3M3 &15 & 0.37 & 984 & 0.76 & 8.58 & 9.34 & 2.00 & 5047 & 3.83 &  84.40& 88.23 \\
         &10 & 0.26 &1087 & 0.31 & 8.95 & 9.26 & & & & & \\ \hline

MSC04515 &13 & 0.65 & 712 & 0.28 & 4.16 & 4.44 & 0.66 & 1020 & 0.20 & 5.93 & 6.13\\
         &10 & 0.52 & 809 & 0.22 & 4.47 & 4.69 & & & & & \\ \hline

EX15     &20 & 1.11 & 601 & 0.81 & 4.88 & 5.69 & 1.83 & 1325 & 0.55 & 12.88 & 13.43\\
         &10 & 0.65 & 511 & 0.34 & 3.53 & 3.87 & & & & & \\ \hline

KUU      &8  & 0.19 & 141 & 0.38 & 1.70 & 2.08 & 0.18 & 144 & 0.25 & 1.73 & 1.98\\
         &10 & 0.23 & 131 & 0.42 & 1.67 & 2.09 & & & & & \\ \hline
\end{tabular}
}
\end{table}

In Table 3, the numerical results for matrices BCSSTK13 and
BCSSTK27 with different values of $lfil$ are given. This table
shows the effect of an increase in $lfil$ on the reduction of the
number of the iterations for convergence. The results of this
table also indicate that the choices $eps=0.01$ and $lfil=10$
lead to good results.

%************Table 3   *******************
\begin{table}
\caption{ Results for matrices BCSSTK13 and BCSSTK21 with
different values of  $lfil$ }\vspace{0.4cm}

\centering
\begin{tabular}{|c||c|c|c||c|c|c|}

\hline

& \multicolumn{3}{|c||}{BCSSTK13} &
\multicolumn{3}{|c|}{BCSSTK21}  \\
\hline

$lfil$ &  P-time & It-time & P-Its  &P-time & It-time & P-Its \\
\hline \hline

2  &  0.03 & 3.02 & 1039 & 0.06  & 0.81 & 316 \\\hline

4  & 0.05  & 2.77 & 917  & 0.08  & 0.75 & 270 \\\hline

6  & 0.06  & 2.48 & 793  & 0.09  & 0.72 & 246  \\\hline

8  & 0.08  & 2.19 & 685  & 0.10  & 0.61 & 198  \\\hline

10 & 0.11  & 1.81 & 550  & 0.11  & 0.53 & 164  \\\hline

12 & 0.14  & 1.69 & 514  & 0.14  & 0.50 & 148  \\\hline

14 & 0.17  & 1.84 & 529  & 0.17  & 0.50 & 142  \\\hline

16 & 0.22  & 1.78 & 502  & 0.20  & 0.47 & 125  \\\hline

\end{tabular}
\end{table}

In \cite{SAINV}, the numerical results of the SAINV preconditioner
in conjunction with some preliminary transformations operated on
the coefficient matrix such as symmetric diagonal scaling,
reordering with the multiple minimum degree (MMD) algorithm
\cite{Liu} and diagonally compensated reduction of positive
off-diagonal entries (DCR), were given. The authors have
concluded that the SAINV preconditioner in conjunction with
symmetric diagonal scaling and reordering with MMD (J-MMD-SAINV)
is often the best choice between the variants of the
preconditioners used in \cite{SAINV}. In continuation, we use 14 out
of 16 matrices used in \cite{SAINV} for the numerical
experiments. Most of these matrices can be extracted from the Matrix
Market website. The exceptions are NASA2910 and NASA4704 which
can be downloaded from Tim Davis's collection, and the SMT
matrix, which was provided by R. Kouhia of the Helsinki University
of Technology
\footnote{http://users.tkk.fi/$\sim$kouhia/sparse.html}. In Table
4, we give the numerical results of the new preconditioner in
conjunction with symmetric diagonal scaling (J-N-M) and the
J-MMD-SAINV preconditioner. Results of the J-MMD-SAINV
preconditioner and the number of iterations of the conjugate
gradient algorithm with Jacobi preconditioning (JCG-Its) were
extracted from Table 7 and Table 1 in \cite{SAINV}, respectively.
It is necessary to mention that the parameter $lfil$ was chosen
such that the parameter $\rho$ of the new preconditioner is less
than or approximately equal to that of the SAINV preconditioner.
All assumptions and notations are as before.

Table 4 shows that the new preconditioner is better than the
J-MMD-SAINV preconditioner for 9 out of 14 matrices presented in
this table.

%************Table 4  *******************
\begin{table}
\caption{Numerical results of the new method in conjunction with
symmetric diagonal scaling and the  J-MMD-SAINV
preconditioner.}\vspace{0.4cm}

\centering {\tiny
\begin{tabular}{|l||c|c||c|c|c|c|c|c||c|c||c|  }
\hline &\multicolumn{2}{|c||}{}&\multicolumn{6}{|c||}{J-N-M} &
\multicolumn{2}{|c||}{J-MMD-SAINV}& \\
\hline \hline

~~~~matrix & $n$ & $nnz$ & $lfil$&$\rho$ & P-Its & P-time & It-time&
T-time &$\rho$& P-Its & JCG-Its   \\
\hline \hline

BCSSTK13 & 2003 & 42943 & 17 & 0.38 & 275 & 0.17 & 0.94 & 1.11 & 0.39 & 349& 1406 \\
\hline

BCSSTK14 & 1806 & 32630 & 9 & 0.28 & 83 & 0.06 & 0.23 & 0.29 & 0.27  &  73 & 409 \\
\hline

BCSSTK15 & 3948 & 60882 & 11 & 0.32  & 176 & 0.17 & 0.95 & 1.12 & 0.33 & 167 & 518\\
\hline

BCSSTK16 & 4884 & 147631 & 6 & 0.12 & 95 & 0.20 & 0.89 & 1.09 &  0.12 &  98  & 191\\
\hline

BCSSTK17 & 10974& 219812 & 16 & 0.40 & 653 & 1.19 & 12.03 & 13.22 & 0.40 & 711 & 2522 \\
\hline

BCSSTK18 & 11948& 80519 & 8 & 0.57 & 515 & 0.78 & 5.73 & 6.51 & 0.58 & 261 & 1120  \\
\hline

BCSSTK21 & 3600 & 15100 &10 & 1.46 & 179 & 0.11 & 0.58 & 0.69 & 1.51 & 191 & 559 \\
\hline

BCSSTK25 & 15439 & 133840 &10 & 0.59 & 1614 & 1.45 & 26.86 & 28.31 & 0.57 & 1512&$\dag$ \\
\hline

S1RMQ4M1 & 5489 & 143300 &8  & 0.19 & 247 & 0.27 & 2.41 & 2.68 & 0.20 & 248 &692\\
\hline

S2RMQ4M1 & 5489 & 143300 &10 & 0.23 & 403 & 0.33 & 3.98 & 4.31 & 0.25 & 528 & 1529\\
\hline

S3RMQ4M1 & 5489 & 143300 &13 & 0.24 & 569 & 0.36 & 5.66 & 6.02 & 0.24 & 1140 &6884\\
\hline

NASA2910 & 2910 & 88603 &20 & 0.32 & 262 & 0.41 & 1.64 & 2.05 & 0.95 & 341 & 1350\\
\hline

NASA4704 & 4704 & 54730 &20 & 0.85 & 567 & 0.47 & 3.77 & 4.24 & 0.91 & 1176 & 4866\\
\hline

SMT      & 25710 & 1889447&15 & 0.11 & 734 & 6.28 & 74.83 & 81.11 & 0.11 &  546 &1984\\
\hline
\end{tabular}
}
\end{table}

The last section of numerical experiments is devoted to some large
matrices extracted from Tim Davis's collection. Numerical results
with different values of $lfil$ and with $eps=0.01$ are given in Table
5. As can be seen, the new preconditioner in conjunction with
symmetric diagonal scaling furnishes  good results for large
matrices.

%************Table 4  *******************
\begin{table}
\caption{Numerical results of the new method in conjunction with
symmetric diagonal scaling.}\vspace{0.4cm}

\centering {\tiny
\begin{tabular}{|l||c|c||c|c|c|c|c|c||c|c|  }
\hline &\multicolumn{2}{|c||}{}&\multicolumn{6}{|c||}{New method with symmetric diagonal scaling}
& \multicolumn{2}{|c|}{JCG} \\
\hline \hline

~~~~matrix & $n$ & $nnz$ & $lfil$&$\rho$ & P-Its & P-time & It-time&
T-time & Its  & Time  \\
\hline \hline

GRIDGENA & 48962 & 280523   & 10 & 1.07 & 609 & 11.36 & 29.85 & 41.21 & 1720  &  48.61\\
                      &  &  & 15 & 1.52 & 540 & 11.98 & 29.34 & 41.32 &       & \\
                      &  &  & 20 & 2.00 & 464 & 13.92 & 27.89 & 41.81 &       & \\
\hline

CVXBQP1 & 50000  & 199984   & 10 & 1.14 & 1053 & 11.38 & 45.33 & 56.71 & 3330  & 90.95\\
                      &  &  & 15 & 1.50 &  886 & 11.80 & 41.00 & 52.80 &       & \\
                      &  &  & 20 & 1.79 &  750 & 12.28 & 36.84 & 49.12 &       & \\
\hline

APACHE1 & 80800  & 311492   & 10 & 1.43 & 325 & 29.44 & 24.20 & 53.64 & 1796 & 79.56\\
                      &  &  & 15 & 1.80 & 245 & 30.36 & 19.48 & 49.84 &      & \\
                      &  &  & 20 & 2.22 & 218 & 31.44 & 18.52 & 49.96 &      & \\
\hline

CF2D2   & 123440 & 1605669  & 10 & 0.45 & 1002 & 68.63 & 200.16 & 268.79 & 3870  & 368.91 \\
                      &  &  & 15 & 0.65 &  850 & 71.17 & 179.8  & 250.97 &   & \\
                      &  &  & 20 & 0.85 &  717 & 73.88 & 164.23 & 238.11 &   & \\
\hline
\end{tabular}
}
\end{table}

\indent We end this section by giving the numerical results for the
matrix APACHE2 extracted from Tim Davis's collection. This is a
large matrix of dimension $n=715176$ with $nnz=2 776 523$ nonzero
entries in the lower part. The conjugate gradient method in
conjunction with symmetric diagonal scaling converges in 2482
iterations. The conjugate gradient method in conjunction with
the new preconditioner and symmetric diagonal scaling with
$lfil=5~(\rho=0.98)$ and $lfil=10~(\rho=1.52)$ converges in $839$
and $575$ iterations, respectively. Convergence history of these
methods is displayed in Figure 1.

Numerical results for the matrix APACHE2 and matrices in Table 5
(and previous tables) show that the new preconditioner reduces the
number of iterations by about a factor of three. This is true for
the J-MMD-SAINV preconditioner based upon a conclusion reported in (
\cite{SAINV}, page 1328).

\begin{figure}
\centering \includegraphics[height=6cm,width=8cm]{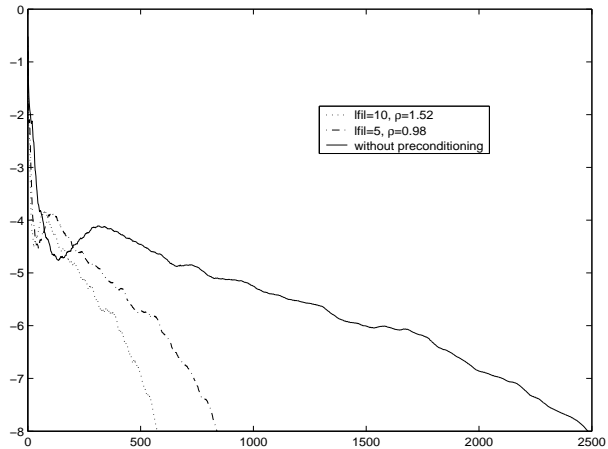}
\caption{Convergence history of the new preconditioner for the
matrix APACHE2.}
\end{figure}

\bigskip

\medskip

\noindent \textbf{ 5. Conclusion and future work}
\bigskip

We have proposed an approach for computing a sparse incomplete
factorization of the inverse of an SPD matrix. The resulting
factorized sparse approximate inverse was used as a preconditioner
for solving symmetric positive definite linear systems of
equations by using the conjugate gradient algorithm.  The new
preconditioner does not need to specify the sparsity pattern of
the inverse factor in advance. For preserving sparsity it is
enough to specify two parameters $eps$ and $lfil$. Numerical
results show that $eps=0.01$ and $lfil=10$ usually give good
results. Our numerical results also show that the proposed method
in conjunction with the symmetric diagonal scaling is somewhat
better than  the J-MMD-SAINV preconditioner.

The new preconditioner is suitable for parallel computers.  It can also
be implemented for  normal equations with a little revision.

Future work may focus on extending the proposed preconditioner  to
general matrices and studying the effect of different reordering
techniques on the convergence rate.

\bigskip

\bigskip

\noindent \textbf{ 6. Acknowledgments}

\bigskip

The authors are gratefully indebted to Edmond Chow for carefully
reading of an earlier draft of this paper and giving several
valuable comments. 
%The authors also are grateful to the anonymous
%referees and our editor Dr. Eric de Sturler for their comments which
%substantially improved the quality of this paper.

\end{document}